\documentclass[aps,prl,preprint,groupedaddress]{revtex4-1}

\usepackage{amsmath,amssymb,amsfonts}
\usepackage{graphicx}
\usepackage{xcolor}
\newcommand{\referee}[1]{\textcolor{red}{#1}}

\begin{document}

\title{Locally Optimal Control of Complex Networks}

\author{Isaac Klickstein}
\email{iklick@unm.edu}
\affiliation{Department of Mechanical Engineering, University of New Mexico}

\author{Afroza Shirin}
\affiliation{Department of Mechanical Engineering, University of New Mexico}

\author{Francesco Sorrentino}
\affiliation{Department of Mechanical Engineering, University of New Mexico}

\date{\today}

\begin{abstract}
It has recently been shown that the minimum energy solution of the control problem for a linear system produces a control trajectory that is nonlocal.
An issue then arises when the dynamics represents a linearization of the underlying nonlinear dynamics of the system where the linearization is only valid in a \textit{local} region of the state space.
Here we provide a solution to the problem of optimally controlling a linearized system by deriving a time-varying set that represents all possible control trajectories parameterized by time and energy.
As long as the control action terminus is defined within this set, the control trajectory is guaranteed to be local.
If the desired terminus of the control action is far from the initial state, a series of local control actions can be performed in series, re-linearizing the dynamics at each new position.
\end{abstract}

\pacs{}

\maketitle
\referee{
Recent work investigates control strategies for complex networks governed by nonlinear dynamical equations.
}
\cite{liu2016control, gao2016universal, ruths2014control, yan2015spectrum, iudice2015structural}.
%
Such problems occur in opinion dynamics in a population \cite{proskurnikov2016opinion}, consensus in robotic networks \cite{wang2014distributed}, developing gene therapies \cite{zhou2014human}, avoiding cascade failures in power grids \cite{pagani2013power}, and many others \cite{pecora2014cluster, sorrentino2016complete, gao2016universal, yan2012controlling, ji2013cluster}.\\
\indent
\referee{
Recently \cite{sun2013controllability}, it was shown that the minimum energy state trajectory of linear systems is \emph{nonlocal} and so one should not attempt to apply minimum energy control to linearized systems.
We instead focus on determining the region of state space where the trajectory does remain local and so minimum energy control can still be applied to linearized approximations of nonlinear systems.
We apply our results to develop an algorithm that determines a piecewise open-loop control signal for nonlinear systems.
}\\
\indent
\referee{The control of complex networks goverened by nonlinear dynamical equations is still in its infancy \cite{liu2016control}.
We attempt to bridge the gap between the well developed techniques for controlling complex networks governed by linear dynamics and those networks which are goverened by nonlinear dynamics.
A recent paper \cite{cornelius2013realistic} suggested that by perturbing the initial state of the system one may be able to place the perturbed initial state in the basin of attraction of a desirable attractor.
A similar method proposed \cite{ott1990controlling, lai2014controlling} applies perturbations to the system parameters rather than to the states like in the previous method.
If the dynamical equations are diffusive \cite{mochizuki2013dynamics, fiedler2013dynamics, zanudo2016structure}, then by over-riding the dynamics of the nodes in a feedback vertex set, one can drive the remaining nodes to an attractor.
State and parameter perturbations represent heuristic methods that attempt to alter either an initial condition or the dynamical equations themselves to move the state into a desirable attractor's basin of attraction.}
\\
\indent
\referee{This letter describes two main results.
The first result is a derivation of a time-varying ellipsoid where all minimum energy state trajectories remain local.
The second result applies the time-varying ellipsoid equation to develop a piecewise controller that drives a nonlinear dynamical network's states towards the basin of attraction of a desired attractor.}\\
%
\indent
\referee{%
It was recently shown \cite{sun2013controllability} that choosing arbitrary initial and final conditions of the minimum energy optimal control problem leads to a \emph{non-local} state trajectory.
Specifically, this means that the length of the state trajectory is independent of the distance between initial and final conditions in average.
Our first result determines the particular set of final conditions that guarantee the locality of the minimum energy controlled state trajectory of a linear dynamical system.
The minimum energy control signal is found by solving the optimal control problem,
\begin{equation}\label{eq:optcon}
  \begin{aligned}
    \min\limits_{\textbf{u}(t)} && &\frac{1}{2}\int_{t_0}^{t_f} \textbf{u}^T(t) \textbf{u}(t)dt\\
    \text{s.t.} && &\dot{\textbf{x}}(t) = A \textbf{x}(t) + \textbf{f} + B \textbf{u}(t)\\
    && &\textbf{x}(t_0) = \textbf{x}_0, \text{ } \textbf{x}(t_f) = \textbf{x}_f
  \end{aligned}
\end{equation}
After computing the Hamiltonian and solving the resulting system of ODEs (see SI section 1), the minimum energy control signal is found to be $\textbf{u}(t) = B^T e^{A^T(t_f-t)} W^{-1}(t_f) (\textbf{x}_f - \textbf{g}_f)$ where $W(t)=e^{A(t-\tau)}BB^Te^{A^T(t-\tau)}d\tau$ is the \emph{controllability Gramian} and $\textbf{g}(t) = e^{A(t-t_0)} \textbf{x}_0 + \int_{t_0}^t e^{A(t-\tau)} d\tau \textbf{f}$ is the \emph{zero-input state trajectory}.
The energy (or effort) consumed by the control signal, $E(t)$, is a monotonically increasing, positive definite function defined as the cumulative sum of squares of each individual signal,
\begin{equation}\label{eq:energy}
  E(t) = \int_{t_0}^t \textbf{u}^T(t) \textbf{u}(t) = \left(\textbf{x}(t) - \textbf{g}(t)\right)^T W^{-1}(t) \left(\textbf{x}(t) - \textbf{g}(t)\right)
\end{equation}
The energy consumed can be expressed as an equation of an $n$-dimensional hyper-ellipsoid centered at $\textbf{g}(t)$ and with principal axes in the eigen-directions of $W(t)$, each with width equal to $2 \sqrt{E(t)} \sqrt{d_i(t)}$ where $d_i(t)$ is the corresponding eigenvalue of $W(t)$.
The hyper-ellipsoid, defined as $\mathcal{S}(t) = \left\{\textbf{x}(t) | \left(\textbf{x}(t) - \textbf{g}(t)\right)^T W^{-1}(t) \left(\textbf{x}(t) - \textbf{g}(t)\right) = E(t) \right\}$ represents the set of states reachable with $E(t)$.
Note that the set of states corresponding to a particular value $E(t_1)$ is \emph{independent} of the function form $E(t)$, $t \in [t_0,t_1)$ for all previous time.
By restricting the amount of energy available, $E(t_f)$, we can determine a set of final conditions, $\mathcal{S}(t_f)$ such that the state trajectory remains \emph{local}, that is, the state at time $t$ lies on the hyper-ellipsoid $\mathcal{S}(t)$.\\
\begin{figure}[t!]
  \centering
  \includegraphics[scale=1]{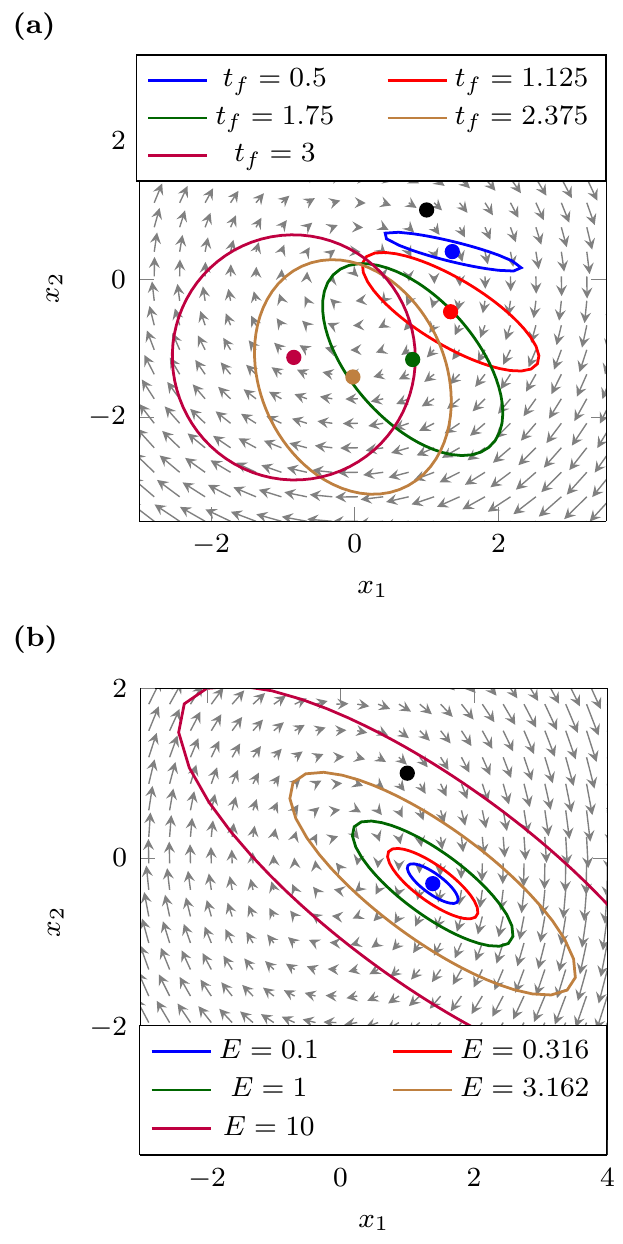}
  \caption{
    Visualizing the ellipsoid by varying $t_f$ and $E$.
    The dynamics are linear, two dimensional with equations $\dot{x}_1=x_2 + u$ and $\dot{x}_2 = -x_1$.
    The arrows represent the flow of the zero-input system.
    In panel (a), the final time is increased from 0.5 to 3.0 while $E$, the energy, is held constant.
    Note how the center of the ellipsoid moves with the zero-input trajectory, i.e., it follows the flow.
    Also, as the final time increases, the direction and width of axes change.
    In panel (b), the final time is held constant while $E$, the energy, is increased.
    The center of the ellipsoid and the directions of the axes remain fixed as $E$ increases and the only aspect that is altered is the width of the axes.
    }
    \label{fig:ellipsoids}
\end{figure}
\indent
To further understand the structure and evolution of the hyper-ellipsoid $\mathcal{S}(t)$, in Fig. \ref{fig:ellipsoids} we consider a linear system consisting of two states, $\dot{x}_1 = x_2 + u$ and $\dot{x}_2 = -x_1$ with initial conditions $x_1(0)=1$ and $x_2(0)=1$.
In Fig. \ref{fig:ellipsoids}(a), the energy $E(t)=1$ is held constant while the time $t$ is allowed to grow.
We see the centroid of the ellipsoid $\mathcal{S}(t)$ moves with the zero-input state trajectory while the axes grow and rotate.
Note that for short $t_f$, the larger axis is primarily in the direction of the state of the node that receives the control input, $x_1$.
On the other hand, in Fig. \ref{fig:ellipsoids}, holding $t$ constant and increasing $E(t)$ does not vary the centroid of the hyper-ellipsoid or the relative axis widths, but only scales the ellipsoid.\\
\indent
We have shown how the hyper-ellipsoid $\mathcal{S}(t_f)$ can be designed by choosing $t_f$ and $E(t_f)$, which in turn allows us to choose a final condition $\textbf{x}(t_f) \in \mathcal{S}(t_f)$ to ensure a local minimum state trajectory.
We will now show how this result can be applied when developing a piecewise controller for a large class of nonlinear dynamical systems.
}\\
\indent
We focus on affine systems of nonlinear differential equations.
The differential equation that describes the behavior of a single node is,
\begin{equation}\label{eq:dyn}
  \dot{x}_i(t) = F_i(\textbf{x}(t)) + \sum_{k=1}^M b_{ik} u_k(t)
\end{equation}
where $i = 1,\ldots,N$ and all functions, $F_i:\mathbb{R}^N \rightarrow \mathbb{R}$ are assumed to be \emph{smooth}.
The sum $\sum_{k=1}^M b_{ik} = 1$ if node $i$ is a \textit{driver node} and $\sum_{k=1}^M b_{ik} = 0$ if node $i$ is not a driver node.
To reflect the network nature of our problem, the values $b_{ik}$, $i = 1,\ldots,N$, $k = 1,\ldots,M$ are either 0 or 1 where if $b_{ik} = 1$ then node $i$ receives input $k$ and if $b_{ik} = 0$, then node $i$ does not receive input $k$, \referee{and each input is received by only one node}.
Equation \eqref{eq:dyn} can be rewritten in vector form, $\dot{\textbf{x}}(t) = \textbf{F}(\textbf{x}(t)) + B \textbf{u}(t)$ where $\textbf{x}(t)$ is the $N$-vector of states for the nodes of the network, $B$ is the $N \times M$ matrix with elements $b_{ik}$ and $\textbf{u}(t)$ is the $M$-vector of control inputs.\\
\indent
\referee{First, we will examine the local problem before demonstrating how a series of local problems can be constructed to develop a piecewise control strategy to drive the network's states into the basin of attraction of a desired attractor.}
The dynamics of a nonlinear system can be approximated locally about a \referee{non-equilibrium} point $\textbf{x}_p$ by a first order Taylor expansion,
\begin{equation}\label{eq:linear}
  \dot{\textbf{x}}(t) = \textbf{f}_p + A_p \textbf{x}(t) + B \textbf{u}(t) + H.O.T.
\end{equation}
where $A_p = \left.\frac{\partial \textbf{f}}{\partial \textbf{x}} \right|_{\textbf{x}=\textbf{x}_p}$ is the Jacobian of the nonlinear dynamics evaluated at $\textbf{x}_p$ and $\textbf{f}_p = \textbf{f}(\textbf{x}_p) - A_p \textbf{x}_p$ represents the flow at $\textbf{x}_p$.
The higher order terms are collected in $H.O.T.$.
This linearization is a valid representation of the nonlinear dynamics only in a local region of state space centered at $\textbf{x}_p$.
We qualify the region where the linearization in Eq. \eqref{eq:linear} is valid as the compact, but not necessarily convex, set,
\begin{equation}\label{eq:VLR}
  \mathcal{N}_p = \left\{ \textbf{x} \in \mathbb{R}^N | \text{ }||\textbf{f}(\textbf{x})- \textbf{f}_p - A_p \textbf{x}||_2 \leq \epsilon \right\}
\end{equation}
where $\epsilon > 0$ is a `small' positive scalar that represents the desired quality of the linearized region, that is, what is the largest deviation between the nonlinear and linear dynamics we may allow.
\referee{It is important to make explicit that the linearization is also temporal.
As a trivial example, assume $\textbf{u}(t) \equiv \boldsymbol{0}$ and set $\textbf{x}(0) = \textbf{x}_p$.
The linearized dynamical equations at this point, $\dot{\textbf{x}}(t) = \textbf{f}(\textbf{x}_p) \neq \boldsymbol{0}$ so at some time $t > 0$, $\textbf{x}(t) \notin \mathcal{N}_p$.
Note that the resulting linear system in Eq. \eqref{eq:linear} appears in the optimal control problem in Eq. \eqref{eq:optcon}.}\\
\indent
\referee{
We have now defined how one can compute a local minimum energy control signal of the linearized model of the true nonlinear system of differential equations such that the state remains within the valid linearization region in Eq. \eqref{eq:VLR} at all times $t \in [t_0,t_f]$.
}
%
\indent We now reframe the above local problem in terms of a global problem.
Let the initial state of the system be denoted $\textbf{x}(0) = \textbf{x}_0$ and some desired region of state space $\mathcal{X} \in \mathbb{R}^N$ which we want to reach in finite time.
A typical example defines $\mathcal{X}$ as a conservative approximation of the basin of attraction of a desirable attractor of $\textbf{f}(\textbf{x}(t))$.
After computing the linearized dynamics at $\textbf{x}_0$, that is $\textbf{f}_0$ and $A_0$, one may choose a time $t_1$ and a point $\textbf{x}_1 \in \mathcal{S}_0(t_1) \in \mathcal{N}_0$ so that the control trajectory is entirely inside the valid linearization neighborhood, Eq. \eqref{eq:VLR}, using the methods described above (and elaborated upon in the SI).
We can then re-linearize the system about $\textbf{x}_1$, that is compute $\textbf{f}_1$ and $A_1$, set the initial time to be $t_1$ and initial state $\textbf{x}_1$, and choose a final time $t_2$ and point $\textbf{x}_2 \in \mathcal{S}_1(t_2) \in \mathcal{N}_1$.
This process may be repeated until $\textbf{x}_P(t_P) \in \mathcal{X}$ at iteration $P$.
The returned solution is a series of times $t_p$ and points $\textbf{x}_p$, $p = 0,1,\ldots,P$.\\
\indent Two important caveats must be stated with respect to the previous iterative approach.
The first is that there is no guarantee that a series of points $t_p$ and $\textbf{x}_p$ such that $\textbf{x}_p \in \mathcal{S}_{p-1}(t_p) \in \mathcal{N}_{p-1}$, $p = 1, \ldots,P$ and $\textbf{x}_P \in \mathcal{X}$ exists.
The method presented above is only a guarantee that a minimum energy control trajectory remains local if the terminal point is chosen in the ellipsoid $\mathcal{S}(t_f)$ that is completely contained in the neighborhood $\mathcal{N}$.
The second caveat is that the choice of each point, $\textbf{x}_p$, may or may not allow for achieving the desired final condition.
The decision mechanism to choose each next point $\textbf{x}_{p+1}$ we have used in practice collects $Q$ feasible choices, $\textbf{x}_{p+1}^{(k)}$, $t_{p+1}^{(k)}$, with their respective energy requirement, $E^{(k)}$, $k = 1,\ldots,Q$, and assigns a \textit{fitness} to each one,
\begin{equation}\label{eq:fitness}
  F^{(k)} = k_E E^{(k)} + \sum_{i=1}^N w_i \left((\textbf{x}^{(k)}_{p+1})_i - \left(\textbf{x}_f\right)_i\right),
\end{equation}
where $k_E$ is a positive weight applied to the required energy for the choice and $w_i$, $i = 1,\ldots,n$, is a weight applied to the remaining change required for each node's state.
The next point of the iterative process, $\textbf{x}_{p+1} = \textbf{x}_{p+1}^{(k)}$, is chosen such that $ k = \text{argmin}_k F^{(k)}$.\\
\begin{figure}
  \centering
  \includegraphics[scale=1]{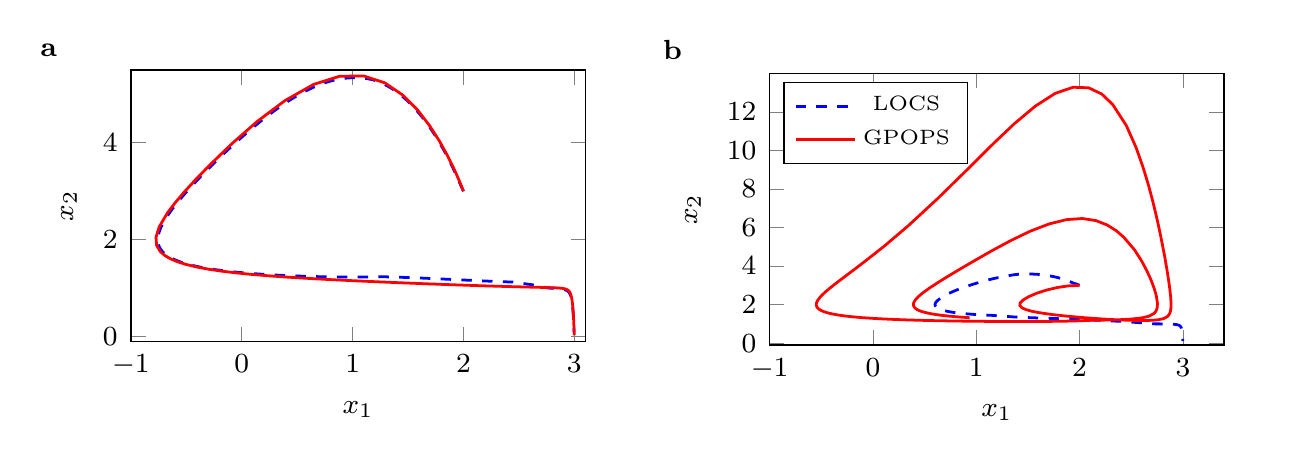}
  \caption{A two-dimensional example, with exact knowledge of the system and with imperfect knowledge of the system.
  (a) State trajectory when the system is known exactly
  The fitness function is set to $F = 4E + 0.8|x_1 - x_{f,1}| + 0.2 |x_2 - x_{f,2}|$.
  To compare, the solution returned from GPOPS is also shown which, as can be seen, is similar to the routine returned from LOCS.
  Note that this is a result of choosing the fitness function carefully based on the dynamics of the system.
  (b) State trajectory when the system is not known exactly.
  In this case GPOPS, which gives a trajectory planning solution, is inadequate (not robust).
  On the other hand, as LOCS re-evaluates the linearization at each step, robustness is built into the method and can compensate for the model uncertainty.
  }
  \label{fig:comb}
\end{figure}
\indent To demonstrate an implementation of LOCS, we first consider a two-dimensional system to enhance the visualization.
Further details of the implementation are described in the Supplementary Information.
The two-dimensional nonlinear system is governed by nonlinear differential equations $\dot{x}_1 = (x_1-3)(x_2-2)$ and $\dot{x}_2 = x_2(x_1-1)(x_2-1) + u$.
The global control action tries to move the system from a periodic orbit about the fixed point at (1,2) to a positively invariant set around the stable fixed point at (3,0).
In order to verify the solution at each step, the minimum energy control input is applied to the nonlinear dynamical equations using a Runge-Kutta (RK) solver to simulate the dynamics.
The state at the end of each step returned by the RK solver, $\textbf{x}(t_p)$, is used as the subsequent initial condition for step $(p+1)$.
Using this method, the small imperfections of the linear model are retained.
Figure \ref{fig:comb}(a) shows the returned state trajectory from the algorithm as a blue dashed curve.
Also included is the minimum energy trajectory returned from GPOPS (red curve).
GPOPS \cite{patterson2014gpops} is a well-known direct method to solve optimal control problems by discretizing the states and inputs in time and solving the resulting nonlinear programming problem.
With knowledge of the dynamics, we choose the \emph{fitness function} for each possible point $\textbf{x}(t_{p+1}) = \textbf{x}_{p+1}$ to be $F = 4E + 0.8|x_1-x_{f1}| + 0.2|x_2-x_{f2}|$.
At each step, we sample 40 feasible choices for the next step.\\
\indent
Another benefit of the LOCS algorithm is that it is able to correct for the model discrepancies of the linear model.
We explore the ability of LOCS to handle larger imperfections beyond the linearization error.
In Fig. \ref{fig:comb}b, one of the parameters in our dynamical model which we use to compute the linearized model is incorrect.
More specifically, $\dot{x}_1 = (x_1-a)(x_2-2)$, where $a = 3$ in the model and $a = 3.5$ in the actual dynamical equations.
The LOCS algorithm is able to adjust for this parameter inconsistency as the system is relinearized about each successive point returned by the RK simulation.
The re-linearization step builds in robustness as the linear model is updated from step to step during the LOCS procedure.
For comparison, while GPOPS provides a near optimal, trajectory planning solution, it fails at producing satisfactory solutions when the model is not known exactly.
The comparison of LOCS and GPOPS with model uncertainty is shown in Fig. \ref{fig:comb}b with the values of $a$ discussed previously.\\
\indent
\referee{
Here and in the supplementary information, We re-examine two examples recently presented in the literature of the control of complex nonlinear networks.
Both examples use methods referred to as `brute-force control' \cite{lai2014controlling}, that is, there is no clear way one may implement the resulting controllers.
Using our method, we determine a time-varying control signal $\textbf{u}(t)$ that performs the same control action for both examples.\\
\begin{figure}
  \centering
  \includegraphics[scale=1]{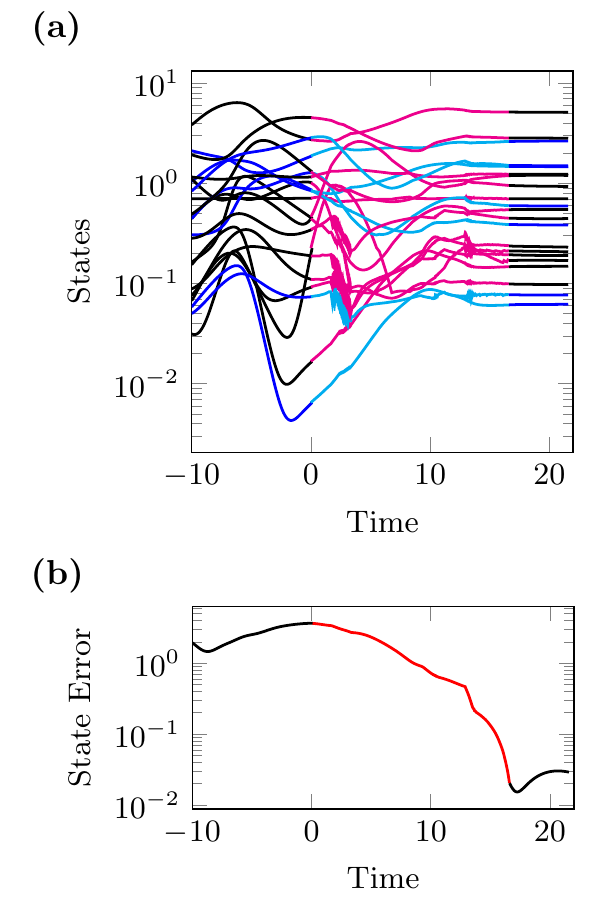}
  \caption{
    Applying the LOCS algorithm to a Mammalian Circadian Rhythm system.
    (a) The state time-evolution.
    The darker colors depict the states during which no control signal is applied to the system.
    The lighter colors represent the states during which the LOCS controller is active.
    Dark blue and light blue curves depict the state evolution of driver nodes while black and magenta curves depict the state evolution of non-driver nodes.
    (b) The norm of the error of the states with respect to the desired fixed point.
  }
  \label{fig:mcr}
\end{figure}
\indent
A model that describes the regulatory structure of the intracellular circadian clock in mice was derived as a system of nonlinear differential equations \cite{mirsky2009model}. 
The system was shown to have a number of attractors, both periodic orbits and fixed points.
This system was used as a demonstration of the utility of the feedback vertex set method \cite{fiedler2013dynamics}.
Rather than over-riding the state variables, we instead attach a control input to each of the states determined to be in the feedback vertex set, which we verified ensures the system is structurally controllable \cite{liu2011controllability}.
With this modification, the complex nonlinear network is in the form required to use LOCS as defined in Eq. \eqref{eq:dyn}.
We define the desired control action as moving from one of the periodic orbits to a stable fixed point.
The state trajectories are plotted in Fig. \ref{fig:mcr}a, where the darker colors correspond to the states before and after the controller is active and the lighter colors correspond to the state while LOCS is active.
The norm of the error as a function of time is shown in Fig. \ref{fig:mcr}b.
}\\
\indent
\referee{
Another relevant example of an application of the LOCS algorithm is to control a network of generators that make-up a power grid after a local failure recently studied in \cite{nishikawa2015comparative,bernstein2014power}.
This is presented in detail in the supplementary information.
}\\
%
\indent In this letter, we presented a methodology to choose a terminal point of a control action for a linearized system such that the optimal control trajectory remains local.
This was accomplished by defining the energy ellipsoid derived from the expression for the control input associated with minimum energy control.
A longer control action can then be defined such that the terminus of one control action becomes the initial state of the next, where the linear dynamics are adjusted by re-linearizing the original nonlinear system about the new initial condition.
Moreover, the algorithm is amenable to a real-time implementation as computations are carried out locally at each point.
For large dimensional systems, such as the dynamical networks we consider in our examples, the LOCS approach provides an open-loop controller that only requires discrete measurements of the states, rather than a continuous feedback loop.
This \textit{daisy chain approach} provides the possibility of using linear optimal control iteratively to traverse larger regions of state space for systems with nonlinear dynamics.
%
\begin{acknowledgments}
The authors thank Franco Garofalo, Francesco Lo Iudice, Karen Blaha and Andrea L'Afflitto for helpful conversations and providing insight into this research.
This work was funded by the National Science Foundation through NSF grant CMMI- 1400193, NSF grant CRISP- 1541148 and ONR Award No. N00014-16-1-2637.
\end{acknowledgments}
\end{document}